\documentclass[graybox]{svmult}

\usepackage{type1cm}        
\usepackage{makeidx}         
\usepackage{graphicx}        
\usepackage{multicol}        
\usepackage[bottom]{footmisc}

\usepackage{tikz,pgfplots} 
\pgfplotsset{compat=newest} 
\pgfplotsset{plot coordinates/math parser=false} 
\usepgfplotslibrary{groupplots}
\usepackage{pgfplotstable}
\newenvironment{customlegend}[1][]{%
	\begingroup
	\csname pgfplots@init@cleared@structures\endcsname
	\pgfplotsset{#1}%
}{  \csname pgfplots@createlegend\endcsname
	\endgroup
}

\def\addlegendimage{\csname pgfplots@addlegendimage\endcsname}

\usepackage{newtxtext}       %
\usepackage{newtxmath}       

\usepackage[software,hardware]{mymacros}
\usepackage{amsbsy} 
\usepackage{bbm}

\newcommand{\bSigma}{\boldsymbol{\Sigma}}
\newcommand{\bhSigma}{\boldsymbol{\hat{\Sigma}}}

\newcommand{\bomega}{\boldsymbol{\omega}}
\newcommand{\bLambda}{\boldsymbol{\Lambda}}
\newcommand{\bOmega}{\boldsymbol{\Omega}}
\newcommand{\bZero}{\mathbf{0}}
\newcommand{\bcM}{\boldsymbol{\cM}}
\newcommand{\bhM}{\mathbf{\hM}}
\newcommand{\bhK}{\mathbf{\hK}}
\newcommand{\bhD}{\mathbf{\hD}}
\newcommand{\bhB}{\mathbf{\hB}}
\newcommand{\bhC}{\mathbf{\hC}}
\newcommand{\bhA}{\mathbf{\hA}}
\newcommand{\bhE}{\mathbf{\hE}}
\newcommand{\bhH}{\mathbf{\hH}}

\newcommand{\bone}{\mathbbm{1}}

\newcommand{\SO}{\texttt{SO}}
\newcommand{\FO}{\texttt{FO}}
\newcommand{\Loew}{\mathbb{L}}
\newcommand{\Loews}{\mathbb{L}_{\sigma}}

\makeindex             

\usepackage{cleveref}
\begin{document}

\title*{Data-Driven Identification of Rayleigh-Damped Second-Order Systems}
\author{Igor Pontes Duff, Pawan Goyal, and Peter Benner}
\institute{Igor Pontes Duff, Pawan Goyal are with the Max Planck Institute for Dynamics of Complex Technical Systems, Sandtorstra\ss e 1, 39106 Magdeburg, Germany; (\email{\{pontes, goyalp\}@mpi-magdeburg.mpg.de}).\\ 
Peter Benner  is with  the Max Planck Institute for Dynamics of Complex Technical Systems, Sandtorstr. 1, 39106 Magdeburg, Germany, and also with the Technische Universit\"at Chemnitz, Faculty of Mathematics, Reichenhainer Stra{\ss}e 41, 09126 Chemnitz, Germany (email: benner@mpi-magdeburg.mpg.de).}
%
%
\maketitle

\abstract{In this paper, we present a data-driven approach to identify second-order systems, having internal Rayleigh damping. This means that the damping matrix is given as a linear combination of the mass and stiffness matrices. These systems typically appear when performing various engineering studies, e.g., vibrational and structural analysis. In an experimental setup, the frequency response of a system can be measured via various approaches, for instance, by measuring the vibrations using an accelerometer. As a consequence, given frequency samples, the identification of the underlying system relies on rational approximation. 
To that aim, we propose an identification of the corresponding second-order system, extending the Loewner framework for this class of systems. The efficiency of the proposed method is demonstrated by means of various numerical benchmarks.
}

\vspace{0.5cm}
\section{Introduction}
In this paper, we discuss a data-driven identification framework for  a class of second-order (\SO) systems of the form:
\begin{equation}\label{eq:SecOrderSys}
\Sigma_{\SO}: = \left\{\begin{array}{r}
 \bM \ddot{\bx}(t) + \bD \dot{\bx}(t) + \bK \bx(t) =  \bB\bu(t), \\  \by(t) = \bC\bx(t), 
\end{array}\right.
\end{equation}   
where $\bx(t) \in \R^{n}$ is the state vector,  $\bu(t)\in \R^m$ are the inputs, $\by(t)\in\R^p$  are the outputs or measurements, and  $\bM, \bD,\bK \in \R^{n \times n}$ are, respectively, the mass matrix, the damping matrix and the stiffness matrix, $\bB \in \R^{n \times m}$ and   $\bC \in \R^{p \times n}$. For simplicity, we address the problem for single-input single-output (SISO) systems, i.e., $m = p  =1$. The multiple-input multiple-output (MIMO) generalization is straightforward and can be done by following the lines of the MIMO extension of the classical Loewner framework \cite{morMayA07} based on tangential interpolation. Such systems arise in many engineering applications, including vibration analysis \cite{Mei97Princ},  structural dynamics \cite{Cra06Fund} and electric circuits.
We denote the \SO~systems \eqref{eq:SecOrderSys} by $\Sigma_{\SO} = (\bM,\bD, \bK,\bB,\bC)$. Moreover, we assume a zero inhomogeneous condition, i.e., $\bx(0) = \dot{\bx}(0) = 0$. Hence, by means of the Laplace transform, the input-output behavior of the system $\Sigma_{\SO}$ is associated with the transfer function as follows:
\vspace{0.2cm}
\begin{equation}\label{eq:SO_TF} 
\bH_{\SO}(s) = \bC\left(s^2\bM +s\bD +\bK\right)^{-1}\bB. 
\vspace{0.2cm}
\end{equation}
Furthermore, throughout the paper, we assume the proportional Rayleigh damping hypothesis, i.e., the damping matrix $\bD$  is given by a linear combination of the mass and stiffness matrices: 
\vspace{0.2cm}
\begin{equation}\label{eq:Rayleigh_damping} 
\bD = \alpha \bM + \beta \bK,
\vspace{0.2cm}
\end{equation}
for $\alpha, \beta \geq 0$. This hypothesis is often considered in several engineering application, where the damper is numerically constructed in order to avoid non-dampened oscillations, see \cite{Mei97Princ} for more details.

In the past twenty years, model order reduction of \SO~systems has been investigated extensively; see for instance  \cite{morMeyS96, morChaLVD06, morReiS08} for balancing-type methods, and \cite{chahlaoui2005model, morBeaG09,wyatt2012issues, Beattie2014h2} for moment matching and $\cH_2$-optimality based methods.  Recently, the authors in \cite{saak2019comparison} provided an extensive comparison among common methods for \SO~model order reduction applied to a large-scale mechanical artificial fishtail model.  In all of the above-mentioned works, the authors suppose that they have access to the matrices, defining the original systems and the reduced-order systems are constructed via Petrov-Galerkin projections. Thus, the main goal is to find projection matrices $\bV, \bW \in \R^{n \times r}$, leading to the \SO~reduced-order system
\vspace{0.2cm}
 \begin{equation}\label{eq:SO_TFROM} 
\bhH_{\SO}(s) = \bhC\left(s^2\bhM +s\bhD +\bhK\right)^{-1}\bhB, 
\end{equation}
\vspace{0.2cm}
with $\bhM = \bW^T\bM\bV, \bhD = \bW^T\bD\bV, \bhK = \bW^T\bK\bV, \bhB = \bW^T\bB$ and $\bhC = \bC\bV$.
\vspace{0.1cm}

However, it is not necessary that the realization is given or is feasible to obtain; thus, we suppose that the original system realization may not be available. Instead, we assume to  have access only to frequency domain data, e.g., arising from experiments or numerical simulations.  More precisely, we are interested in solving the following problem.  \newpage
\begin{svgraybox}
	\vspace{-0.3cm}
\begin{problem}[\SO~data-driven identification]\label{pb:SOGeneral} Given interpolation data
	\begin{equation}\label{eq:interdata}
	\{(\sigma_i, \bomega_i) |~\sigma_i \in \mathbb{C}\,~\, \text{and}\,~\, \bomega_i \in \mathbb{C},~i =1,\dots, \rho \},
	\end{equation}
construct a \SO~realization $\Sigma_{\SO} = (\bM,\bD, \bK,\bB,\bC)$ of appropriate dimensions, satisfying the proportional Rayleigh damping hypothesis, i. e.,  \[\bD = \alpha \bM + \beta \bK,\]
whose transfer function $\bH_{\SO}(s) := \bC(s^2\bM+s\bD+\bK)^{-1}\bB$  satisfies the interpolation conditions, i.e.,
	\begin{equation}
	\bH_{\SO}(\sigma_i)= \bomega_i,~ i = 1,\dots \rho.
	\label{eq:constRight} 
	\end{equation} 
\end{problem} 
\end{svgraybox}

Problem~\ref{pb:SOGeneral} corresponds to an identification problem which aims at determining a \SO~realization that not only interpolates at given measurements, but also satisfies the Rayleigh damping hypothesis. A similar problem for time-delay systems was studied in \cite{Pon15real} and \cite{SchU16data}. Furthermore, we would like to mention that a data-driven approach for structured non-parametric systems has been  studied in \cite{schulze2018data}. However, the construction  of the structured reduced-order system is not a straightforward task.


The purpose of this paper is thus to extend the application domain of the Loewner framework established in \cite{IonitaPhd2013,morMayA07} to \SO~systems. With this aim, a new \SO~Loewner framework is developed, yielding a Rayleigh damped \SO~system of the form \eqref{eq:SO_TF} that interpolates at given frequency measurements.  

The rest of the paper is organized as follows.  Section \ref{sec:ClassicLoewnerFram} recalls some preliminary results on the rational interpolation Loewner framework proposed in \cite{morMayA07}. Section \ref{sec:SOLoewnerFramework}  presents an extension of these results to the  class of Rayleigh damped \SO~systems. The section is divided into two parts. The first one assumes the knowledge of the Rayleigh damping parameters, $\alpha$ and $\beta$, and derives the Loewner matrices for \SO~systems. The second part presents a heuristic procedure, originally proposed in \cite{SchU16data} in the context of time-delay systems, enabling us to estimate the parameters $\alpha$ and $\beta$. Finally, Section \ref{sec:NumRes} illustrates the proposed framework by numerical examples and Section \ref{sec:Conc} concludes the paper.

\vspace{-0.3cm}
\section{Classical Loewner framework}\label{sec:ClassicLoewnerFram}
In this section, we briefly recall the Loewner framework \cite{morMayA07}. A first-order (\FO) system $\bSigma_{\FO} = (\bE, \bA, \bB, \bC)$  is a dynamical system of the form:
\begin{equation}\label{eq:FirstOrderSys}
\Sigma_{\FO}: = \left\{\begin{array}{rl}
\bE \dot{\bx}(t)  &= \bA \bx(t) +  \bB\bu(t), ~\,\bx(0) = 0, \\  \by(t) &= \bC\bx(t),
\end{array}\right.
\end{equation}   
with  $\bE, \bA \in \R^{n \times n}$,  $\bB \in \R^{n \times m}$ and $\bC \in \R^{p \times n}$, and the leading dimension $n$ is the order of the system. For clarity of exposition, we focus for now on the single-input single-output (SISO) case, i.e., when $m = p = 1$. The system~\eqref{eq:FirstOrderSys} is associated with the transfer function given by \vspace{-0.15cm}
\begin{equation}\label{eq:FirstOrderTF}
\bH_{\FO}(s) =  \bC\left(s\bE - \bA\right)^{-1}\bB.
\vspace{-0.15cm}
\end{equation}

There exist several MOR techniques for first-order systems such as explicit moment matching \cite{yousuff1985linear,morVilS87}, implicit moment matching using Krylov subspaces \cite{morGalGV96, morGri97}, Sylvester equations based method \cite{morGalVV04a}, extensions for MIMO systems \cite{morGalVV04}. We refer the reader to the books \cite{morAnt05,morBenCOetal17}  for more details. However, our goal lies in the identification of linear systems using only the frequency data. Hence, the identification problem, in its SISO form, is stated as follows.
\begin{svgraybox}
		\vspace{-0.3cm}
\begin{problem}[First-order data-driven model reduction] \label{prob:FOinter} Given interpolation data
	\begin{equation}\label{eq:interdata}
	\{(\sigma_i, \bomega_i) |~\sigma_i \in \mathbb{C}\,~\, \text{and}\,~\, \bomega_i \in \mathbb{C},~i =1,\dots, \rho \} 
	\end{equation}
	construct a minimal-order realization $\bSigma = (\bE, \bA, \bB, \bC)$ of appropriate dimensions, whose transfer function $\bH_{\FO}(s) = \bC(s\bE-\bA)^{-1}\bB$  satisfies the interpolation conditions
	\begin{equation}\label{eq:SISO_FO_inter} 
	\bH_{\FO}(\sigma_i)= \bomega_i,~ i = 1,\dots \rho.
	\end{equation} 
\end{problem}
	\vspace{-0.3cm}
\end{svgraybox}
A wide range of methods has been developed to solve Problem \ref{prob:FOinter}, e.g., vector fitting  \cite{Gus99vecfitt}, the AAA algorithm \cite{Nak18aaa} and the Loewner framework \cite{morMayA07}. In this paper, we focus on the latter approach and, in what follows, we recall some of the results contained therein. Firstly, we assume that the number of interpolation data is even, i.e., $\rho = 2\ell$, and as a result, the data can be partitioned in two disjoint sets as follows:
\vspace{-0.5cm}
\begin{subequations}
\begin{align}\label{eq:rightinterdata}
&\text{right interpolation set  $\cP_r$:}\{(\lambda_i, \mathbf{w}_i) |~\lambda_i \in \mathbb{C}\,~\,\text{and}\,~\,\mathbf{w}_i,~ i =1,\dots, \ell  \},~\text{and}  \\
&\text{left interpolation set $\cP_l$:}\{(\mu_j, \mathbf{v}_j) |~\mu_j \in \mathbb{C}\,~\, \text{and}\,~\, \mathbf{v}_j \in \mathbb{C},~j =1,\dots, \ell \}.\label{eq:leftinterdata}
\end{align}
\end{subequations}
Using this partition, we associate the following Loewner matrices. 
\begin{svgraybox}
	\vspace{-0.1cm}
\begin{definition}[Loewner matrices  \cite{morMayA07}]  Given the right $\cP_r$ and left $\cP_l$ interpolation sets, we associate them with the Loewner matrix $\Loew$ and shifted Loewner matrix $\Loews$ given by
\begin{equation}\label{eq:ClassicLoew}
\Loew = \left( \begin{array}{c c c}
\frac{\bv_1 - \bw_1}{\mu_1-\lambda_1} & \cdots & \frac{\bv_1 -  \bw_{\ell }}{\mu_1-\lambda_{\ell} } \\
\vdots & \ddots & \vdots\\
\frac{\bv_{\ell}-  \bw_1}{\mu_{\ell}-\lambda_1} & \cdots & \frac{\bv_{\ell} -  \bw_{\ell }}{\mu_{\ell}-\lambda_{\ell}} \\
\end{array} \right), \quad \Loews = \left( \begin{array}{c c c}
\frac{\mu_1 \bv_1 - \lambda_1  \bw_1}{\mu_1-\lambda_1} & \cdots & \frac{ \mu_1 \bv_1 - \lambda_{\ell} \bw_{\ell}}{\mu_1-\lambda_{\ell}} \\
\vdots & \ddots & \vdots\\
\frac{\mu_{\ell}\bv_{\ell}- \lambda_1  \bw_1}{\mu_{\ell}-\lambda_1} & \cdots & \frac{\mu_{\ell} \bv_{\ell}  - \lambda_{\ell}  \bw_{\ell}}{\mu_{\ell}-\lambda_{\ell}} \\
\end{array} \right).
\end{equation}
\end{definition}
\vspace{-0.2cm}	
\end{svgraybox}

\begin{remark} The Loewner matrix $\Loew$ was introduced in \cite{Ant86scalar}. As shown therein, its usefulness derives from the fact that its rank is equal to the order of the minimal realization $\bH_{\FO}$ satisfying the interpolation conditions in \eqref{eq:SISO_FO_inter}. Hence, it reveals the complexity of the reduced-order model solving Problem \ref{prob:FOinter}.
\end{remark}

Next, let  us introduce the following matrices associated with the interpolation problem as follows:
\begin{equation}\label{eq:dataMatrices}
\left\{\begin{array}{l}
\bLambda = \diag{\lambda_1, \dots, \lambda_{\ell}} \in \mathbb{C}^{\ell \times {\ell}} \\
\bhH(\bLambda) = \begin{bmatrix}
\mathbf{w}_1 & \dots & \mathbf{w}_{\ell}
\end{bmatrix}^T \in \mathbb{C}^{\ell \times 1}
\end{array}\right.~\,~\textnormal{and}~\,~
\left\{\begin{array}{l}
\bcM = \diag{\mu_1, \dots, \mu_{\ell}} \in \mathbb{C}^{\ell\times \ell} \\
\bhH(\bcM) = \begin{bmatrix}
\mathbf{v}_1 & \dots & \mathbf{v}_{\ell} 
\end{bmatrix}^T\in \mathbb{C}^{\ell \times 1}
\end{array}\right.
\end{equation}
Also, let  $\bone\in \R^{\ell \times 1}$ be the column vector with all entries equal to one. Hence, the Loewner matrices satisfy the following Sylvester equations
\begin{subequations}\label{eq:SylvFOLoew}
\begin{align}\label{eq:SylvFO}
\bcM\Loew -\Loew\bLambda &= \bhH(\bcM)\bone^T - \bone\bhH(\bLambda)^T,~\, \text{and} \\ 
\bcM\Loews -\Loews\bLambda &= \bcM\bhH(\bcM)\bone^T - \bone\bhH(\bLambda)\bLambda.
\end{align}
\end{subequations}
An elegant solution for Problem \ref{prob:FOinter} based on the Loewner pair $(\Loew, \Loews)$ was proposed in \cite{morMayA07}. This is summarized in the following theorem.
\begin{svgraybox}
	\vspace{-0.3cm}
\begin{theorem}[Loewner framework \cite{morMayA07}] Let $\Loew$ and $\Loews$ be the Loewner matrices associated with the partition in \eqref{eq:dataMatrices}. If  $(\Loews, \Loew)$ is a regular pencil  with no $\mu_i$ or $\lambda_j$ being an eigenvalue, then the matrices
\begin{equation*}
\bhE= -\Loew,\quad \bhA= -\Loews, \quad \bhB=\bV, \quad \bhC=\bW, 
\end{equation*}
provides a realization $\bhSigma_{\FO} = (\bhE,\bhA,\bhB, \bhC)$ for a minimal order interpolant of Problem \ref{prob:FOinter}, i.e., the transfer function
\begin{equation*}
\bhH_{\FO}(s) = \bW (s\Loews-\Loew)^{-1}\bV 
\end{equation*}
satisfies the interpolation conditions in \eqref{eq:SISO_FO_inter}.\label{thm:LoewnerFO}
\end{theorem}
	\vspace{-0.3cm}
\end{svgraybox}
Theorem \ref{thm:LoewnerFO}  allows to obtain a \FO~system $\bhH = (\bhE,\bhA,\bhB,\bhC)$ whose transfer function interpolates right and left data as stated in Problem \ref{prob:FOinter}. However, when more data than necessary are provided, then the hypothesis of Theorem \ref{thm:LoewnerFO} may not be satisfied. Hence, a singular-value decomposition (SVD) based procedure has been proposed in \cite{morMayA07} to find an $\FO$ system interpolating the frequency data.

Next, recall that a \SO~system $\Sigma_{\SO} = (\bM,\bD, \bK,\bB,\bC)$  can be written as a first-order realization, for instance, as follows:
\[\bH_{\SO\_\FO}(s) = \cC(s\cE-\cA)^{-1}\cB, \]
where 
\[\cE = \begin{bmatrix}
\bI & \bZero \\ 
\bZero   & \bM 
\end{bmatrix},\quad \cA = \begin{bmatrix}
\bZero & \bI \\ 
-\bK & -\bD 
\end{bmatrix}, \quad \cB =\begin{bmatrix}
\bZero \\ 
\bB
\end{bmatrix}\quad \text{and}~\, \cC = \begin{bmatrix}
\bC & \bZero
\end{bmatrix}.  \]
As a consequence, the classical Loewner framework presented in Section \ref{sec:ClassicLoewnerFram} can be employed to find a first-order realization. However, the intrinsic \SO~structure will not be preserved in the identified model. But the classical Loewner framework yields an information about the order of a \SO~realization fitting the data, which is outlined in the following remark.

\begin{remark}[\bf Order of \SO~model]  Let us suppose that the frequency data in Problem \ref{prob:FOinter} and let $\Loew$ be a Loewner matrix given in \eqref{eq:ClassicLoew} constructed with this data. Then,  the order of the \SO~system fitting the data equals $\frac{1}{2}\rank{\Loew}$.
\end{remark}

In the following section, we discuss an extension of the Loewner framework for the class of Rayleigh damped \SO~systems. 


\section{Second-order Loewner Framework}\label{sec:SOLoewnerFramework}
This section contains our main contribution, which presents an extension of the Loewner framework to the class of \SO~Rayleigh damped systems~\eqref{eq:SecOrderSys}.  Here, we also assume that the number of interpolation data is even, i.e., $\rho = 2\ell$, and the data is partitioned into two disjoint sets as in \eqref{eq:rightinterdata}  and \eqref{eq:leftinterdata}.  Moreover, the data is organized into the matrices $\bLambda$,  $\bhH(\bLambda)$, $\bcM$, $\bhH(\bcM)$ as in \eqref{eq:dataMatrices}. This section is divided into two parts. In the first one, we assume to have a priori knowledge of the Rayleigh damping parameters $\alpha$ and $\beta$ and we derive the equivalents of the Loewner matrices \eqref{eq:ClassicLoew} and the Theorem~\ref{thm:LoewnerFO}  to the class of \SO~Rayleigh damped systems. The second part is dedicated to proposing a heuristic procedure to estimate the parameters $\alpha$ and $\beta$ using the frequency data available.

\subsection{Second-order Loewner matrices}

In what follows, we assume that Problem \ref{pb:SOGeneral} has a minimal order $r$ solution $\bH_{\SO}^{\star}$, given by
\begin{equation}\label{eq:SOansatz}
\bH_{\SO}^{\star}(s) = \bC^{\star}\left(s^2\bM^{\star} +s\bD^{\star} +\bK^{\star}\right)^{-1}\bB^{\star},
\end{equation}
with $\bD^{\star} = \alpha \bM^{\star} + \beta \bK^{\star}$. Here, we also assume that the coefficients $\alpha$ and $\beta$ from the Rayleigh-Damped hypothesis are known. 
Then, later in this section, we will show how to construct a realization equivalent to $\bH_{\SO}^{\star}(s)$ that only depends on the frequency data. To that aim, let us first recall a result from \cite{morBeaG09} enabling projection-based structured preserving model reduction.
\newpage
\begin{theorem}[Structure preserving \SO~model reduction~\cite{morBeaG09}]\label{theo:StructPreserv}
		Consider the \SO~transfer function $\bH_{\SO}(s)$ as given in \eqref{eq:SO_TF}. For given interpolation points $\lambda_i$ and $\mu_i$, $i\in \{1,\ldots,\ell\}$, let the projection matrices $\bV$ and $\bW$ be as follows:\vspace{-0.2cm}
			\begin{subequations}\label{eq:ProjMatrices}
			\begin{align}
			\bV &= \begin{bmatrix}
			\left(\lambda_1^2\bM +\lambda_1\bD + \bK \right) ^{-1}\bB, & \dots, &  \left(\lambda_{\ell}^2\bM +\lambda_{\ell}\bD + \bK \right) ^{-1}\bB \end{bmatrix}	 
			\\ 
			\bW &= \begin{bmatrix}
			 \left(\mu_1^2\bM +\mu_1\bD + \bK \right) ^{-T}\bC^T, & \dots, &  \left(\mu_{\ell}^2\bM +\mu_{\ell}\bD + \bK \right) ^{-T}\bC^T
			\end{bmatrix}
			\end{align}	
		\end{subequations}
		Hence, the reduced-order model $\hat{\bH}_{\SO}(s)$  constructed by Petrov-Galerkin projection as in \eqref{eq:SO_TFROM} satisfies the interpolation conditions
		\begin{align*}
		\bH_{\SO}(\lambda_i) &= \bhH_{\SO}(\lambda_i) \quad \text{and} \quad \bH_{\SO}(\mu_i) = \bhH_{\SO}(\mu_i),~\, \text{for $i =1,\dots, \ell$.} 
		\end{align*}
\end{theorem}
The above theorem  allows us to construct a \SO~reduced-order model by interpolation. Let us apply this theorem to the \SO~system $\bH_{\SO}^{\star}(s)$ \eqref{eq:SOansatz}. For this, we will construct the matrix $\bV$ using the interpolation points in $\bLambda$, and the matrix $\bW$ using the interpolation points in $\bcM$. As a consequence, $\bV$ and $\bW$  are, respectively, the solutions of the following matrix equations
\begin{subequations}
\begin{align}
\bM^{\star} \bV \bLambda^2 + \bD^{\star} \bV\bLambda + \bK^{\star}\bV  &= \bB^{\star}  \bone^T \label{eq:SylvesterV}, \quad \text{and} \\ 
\bcM^2 \bW^T \bM^{\star}  +  \bcM\bW^T\bD^{\star}  + \bW^T \bcM &=  \bone \bC^{\star}\label{eq:SylvesterW},
\end{align}
\end{subequations} 
Multiplying the equations on the left \eqref{eq:SylvesterV} and \eqref{eq:SylvesterW}  on the left by $\bW^T$ and $\bV^T$, respectively, one obtains 
\begin{align*}
\bW^T\bM^{\star} \bV \bLambda^2 +\bW^T \bD^{\star} \bV\bLambda + \bW^T\bK^{\star}\bV &= \bW^T\bB^{\star}  \bone^T ,\\
\bcM^2  \bW^T\bM^{\star}\bV + \bW\bD^{\star}\bV \bcM + \bW\bK^{\star}\bV &= \bC^{\star}\bV \bone^T.
\end{align*}
If we set 
\begin{subequations}
\begin{align}
\bhM = \bW^T\bM^{\star}\bV, \quad \bhD = \bW^T\bD^{\star}\bV, \quad \bhK = \bW^T\bK^{\star}\bV, \\ \bhB = \bW^T\bB^{\star} = \bhH(\bcM), \quad \text{and}\quad  \bhC = \bC^{\star}\bV = \bhH(\bLambda)^T, 
\end{align}
\end{subequations}
then the \SO~system $\bhH_{\SO} = (\bhM, \bhD, \bhK, \bhB, \bhC)$ is the reduced-order model obtained by Theorem \ref{theo:StructPreserv}, satisfying the interpolation conditions from Problem \ref{pb:SOGeneral}. Hence, we can rewrite the above equations as follows:
\begin{subequations}
\begin{align*}
\bhM  \bLambda^2 + \bhD \bLambda + \bhK &= \bhH(\bcM)  \bone^T,  \\
 \bcM^2 \bhM+ \bcM \bhD + \bhK &=  \bone \bhH(\bLambda)^T.
\end{align*}
\end{subequations}
Moreover, if we apply the Raylegh-Damped hypothesis, i.e.,  $\bhD = \alpha\bhM+ \beta\bhK,$ we obtain

\begin{subequations}
	\begin{align}
	\bhM\left(\bLambda^2 + \alpha \bLambda\right) +  \bhK \left(\beta\bLambda+ \bI\right) &= \bhH(\bcM)  \bone^T, \label{eq:PropSylvesterSO1}  \\
	\left(\bcM^2 + \alpha\bcM\right)\bhM  +\left(\beta \bcM+ \bI\right) \bhK &=  \bone \bhH(\bLambda)^T. \label{eq:PropSylvesterSO2}
	\end{align}
\end{subequations}
Notice that the above equations can be solved for $\bhM$ and $\bhK$. However, in order to have an analytic expression for the matrices of the reduced-order system in a similar way as for the Loewner matrices \eqref{eq:ClassicLoew}, we need to introduce the following change of variables:
\begin{subequations}
\begin{align}\label{eq:ChangVarSOLoewner}
  \Loew^{\SO} &:= -(\bI+\beta\bcM)\bhM(\bI+\beta\bLambda), \quad  \Loews^{\SO} := (\bI+\beta\bcM)\bhK(\bI+\beta\bLambda),   \\
  \bB^{\SO}  &:= (\bI+\beta\bcM)\bhH(\bcM), \quad   \text{and}  \quad  ~\bC^{\SO} :=\bhH(\bLambda)^T(\bI+\beta\bLambda).
\end{align}
\end{subequations}
Notice that the two realizations  \[\bhSigma_{\SO} = (\bhM,\alpha\bhM+\beta\bhK, \bhK,\bhB, \bhC)\quad  \text{and} \quad \bhSigma_{\SO}^{\text{Loew}} = (-\Loew^{\SO}, -\alpha \Loew^{\SO} + \beta \Loews^{\SO}, \Loews^{\SO},\bB^{\SO}, \bC^{\SO})\]
are equivalent, i.e., they represent the same transfer function. Hence, the realization $\bhH_{\SO}^{\text{Loew}}$ also satisfies the interpolation conditions from Problem \ref{pb:SOGeneral}. Additionally, by a simple computation, we obtain that the matrices $\Loew^{\SO}$ and $\Loews^{\SO}$ satisfy the following equations
\begin{align*}
\Loew^{\SO} \cF(\bLambda) + \Loews^{\SO} &= -\cD(\bcM)\bhH(\bcM)  \bone^T,\\
\cF(\bcM) \Loew^{\SO} + \Loews^{\SO} &=  -\bone \bhH(\bLambda)^T\cD(\bLambda),
\end{align*}
where , for a given matrix $\bOmega$, $\cF(\bOmega) := (\bI+\beta \bOmega)^{-1}(\bOmega^2+\alpha \bOmega)$ and $\cD(\bOmega) := (\bI+\beta \bOmega)$.  As a consequence, 
\begin{subequations}\label{eq:so_sylLoew}
\begin{align}
\Loew^{\SO}\cF(\bLambda) - \cF(\bcM)\Loew^{\SO} = \bone \bhH(\bLambda)^T\cD(\bLambda) - \cD(\bcM)\bhH(\bcM)\bone^T, 
\\ 
\Loews^{\SO}\cF(\bLambda) - \cF(\bcM)\Loews^{\SO} = \bone \bhH(\bLambda)^T\cN(\bLambda) - \cN(\bcM)\bhH(\bcM)\bone^T, 
\end{align}
\end{subequations}
where, for a given matrix $\bOmega$, $\cN(\bOmega) := (\bOmega^2+\alpha \bOmega)$. Notice that the Sylvester equations~\eqref{eq:so_sylLoew} are equivalent to \eqref{eq:SylvFO} for the case of \SO~systems.  Hence, using those equation, one can derive analytic expressions of $\Loew^{\SO}$ and $\Loews^{\SO}$.

\begin{svgraybox}
	\vspace{-0.5cm}
\begin{definition}[\SO~Loewner matrices]  Let us suppose $\alpha$ and $\beta$ are known and let 

\[d(s) := 1+s\beta, \quad n(s) := s^2+\alpha s, \quad  \text{and} \quad  f(s) := \frac{n(s)}{d(s)},\]
	be scalar functions. Then the \SO~Loewner matrices, namely, the \SO~\,Loewner matrix $\Loew^{\SO}$ and the shifted Loewner matrix $\Loews^{\SO}$ are given by
	\vspace{0.3cm}
    \newpage
    
	\begin{eqnarray}\label{eq:SOLoew}
		\Loew^{\SO} = \left( \begin{array}{c c c}
			\frac{d(\mu_1)\bv_1 -  d(\lambda_1)\bw_1}{f(\mu_1)-f(\lambda_1)} & 
			\cdots 
			& \frac{d(\mu_1)\bv_1 -  d(\lambda_{\ell})\bw_{\ell}}{f(\mu_1)-f(\lambda_{\ell})} \\
			\vdots & 
			\ddots & 
			\vdots \\
			\frac{d(\mu_{\ell})\bv_{\ell} -  d(\lambda_1)\bw_1}{f(\mu_{\ell})-f(\lambda_1)} & \cdots & \frac{d(\mu_{\ell})\bv_{\ell} -  d(\lambda_{\ell})\bw_{\ell}}{f(\mu_{\ell})-f(\lambda_{\ell})} \\
		\end{array} \right) ,
	\end{eqnarray}
	
	\begin{eqnarray}\label{eq:SOsLoew}
	\Loews^{\SO}= \left( \begin{array}{c c c}
		\frac{n(\mu_1) \bv_1 - n(\lambda_1)  \bw_1}{f(\mu_1)-f(\lambda_1)} & \cdots & \frac{n(\mu_1) \bv_1 - n(\lambda_{\ell}) \bw_{\ell}}{f(\mu_1)-f(\lambda_{\ell})} \\
		\vdots & \ddots & \vdots\\
		\frac{n(\mu_{\ell}) \bv_{\ell}-n(\lambda_1)  \bw_1}{f(\mu_{\ell})-f(\lambda_1)} & \cdots & \frac{n(\mu_{\ell}) \bv_{\ell}  -n(\lambda_{\ell})  \bw_{\ell}}{f(\mu_{\ell})-f(\lambda_{\ell})} \\
	\end{array} \right) .
\end{eqnarray}
\end{definition}
	\vspace{-0.5cm}
\end{svgraybox}
Moreover, by construction 
\[ \Loew^{\SO} =   -(\bI+\beta\bcM)\bhM(\bI+\beta\bLambda) =  -(\bI+\beta\bcM)\bW^T\bM^{\star}\bV(\bI+\beta\bLambda). \]
Thus, the following remark holds.
\begin{remark} If we have sufficient interpolation data, then $\rank{\bV} = \rank{\bW} = r.$  As a consequence,  the rank of the \SO~Loewner matrix $\Loew^{\SO}$ gives us the order of the Rayleigh damped \SO~minimal realization interpolating the points, since
	\[ \rank{\Loew^{\SO}} = \rank{\bW^T\bM^{\star}\bV} = \rank{\bM^{\star}} = \text{order of the minimal \SO~interpolant}. \]
\end{remark}

We are now able to state the analogue result to Theorem \ref{thm:LoewnerFO} for Rayleigh-damped $\SO$ systems.
\begin{svgraybox}
	\vspace{-0.5cm}
	\begin{theorem}[\SO~data-driven identification]\label{theo:SOLoew} Assume that $\mu_i \neq \lambda_j$ for all $i,j = 1, \dots, \ell$. Additionally, suppose that  $(s^2+\alpha s)\Loew^{\SO} +(\beta s+1) \Loews^{\SO}$ is invertible for all $s = \{\lambda_1, \dots, \lambda_{\ell}\}\cup \{\mu_1, \dots, \mu_{\ell}\}$.  
Then
		\[
		\begin{array}{rl}
			\bhM = - \Loew^{\SO}, \quad \bhK = \Loews^{\SO},\quad  \bhB = (\bI+\beta\bLambda)^{-1}\bV^{\SO} \quad  \bhC = \bW^{\SO} (\bI+\beta\bcM),
		\end{array}    \] 
		and $\bhK = \alpha \bhM + \beta \bhK$ satisfy the interpolation conditions from Problem \ref{pb:SOGeneral}. 
	\end{theorem}
\end{svgraybox}

We now consider the case where more data than necessary are provided, which is realistic for applications. In this case, the assumptions of the above theorem are not satisfied; thus, one needs to project onto the column span and the row span of a linear combination of the two Loewner matrices. More precisely, let the following assumption be satisfied:
\begin{equation}\label{eq:NumRank}
  \rank{\begin{bmatrix}
	\Loew^{\SO} & \Loews^{\SO}
	\end{bmatrix}} = \rank{\begin{bmatrix}
	\Loew^{\SO} \\ \Loews^{\SO}
	\end{bmatrix}} = r 
\end{equation}
Then, we consider the compact SVDs 
\begin{equation}\label{eq:shortSVD}
\begin{bmatrix}
\Loew^{\SO} & \Loews^{\SO}
\end{bmatrix} = \bY_{\rho} \Sigma_l \tilde{\bV}^T \quad \text{and} \quad  \begin{bmatrix}
\Loew^{\SO} \\ \Loews^{\SO}
\end{bmatrix} = \tilde{\bW}\Sigma_r \bX_{\rho}^T.
\end{equation}
Using the projection matrices $\bV_{\rho}$ and $\bW_{\rho}$, we are able to remove the redundancy in the data by means of the following result.
\begin{theorem} The \SO~realization $\bhSigma_{\SO} = (\bhM,\bhD,\bhK,\bhB,\bhC)$ of a minimal interpolant of Problem \ref{pb:SOGeneral} is given as:
\begin{subequations}
\begin{align}
\bhM &= -\bY_{\rho}^T\Loew^{\SO}\bX_{\rho}, \quad \bhK = -\bY_{\rho}^T\Loews^{\SO}\bX_{\rho}, \quad \bhD = \alpha \bhM +\beta \bhK, \\ \bhB &= \bY_{\rho}^T\bhB^{Loew} ,\quad \text{and}~\quad \bhC =  \bhC^{Loew}  \bX_{\rho}.
\end{align}	
Depending on whether $r$ in \eqref{eq:NumRank} is the exact or approximate rank, we obtain either an interpolant or an approximate interpolant of the data, respectively.
\end{subequations}	
\end{theorem}
\vspace{-0.5cm}
\subsection{Optimizing parameters}\label{subsec:OptPar}
In the previous section, we have shown how to construct a \SO~realization for given transfer function measurements and a priori knowledge of the parameters $\alpha$ and $\beta$ from the Rayleigh-damped hypothesis. However, there are several cases, where exact values of $\alpha$ and $\beta$ are not known but we rather can have a hint of the range for the parameters, i.e., $\alpha \in \cR_\alpha$ and $\beta \in \cR_\beta$. Therefore, as done for delay systems in \cite{SchU16data}, we also propose a heuristic optimization approach to obtain the parameters $\alpha$ and $\beta$ for \SO~systems, satisfying the Rayleigh-damped hypothesis. For this purpose, we split the data training $\cD_{\text{training}}$ and test set $\cD_{\text{test}}$, e.g., in the ratio $80{:}20$. Hence, we ideally aim at solving the optimization as follows:
\begin{equation}\label{eq:opt_problem}
\min_{\alpha\in \cR_\alpha,~\,\beta\in \cR_\beta} \cJ(\alpha, \beta)
\end{equation}
where 
$$\cJ(\alpha, \beta) := \sum_{(\sigma_k,v_k) \in \cD_{\text{test}} } \left\|\bhH_{\SO}(\sigma_k\alpha,\beta) - v_k\right\|^2 + \sum_{(\mu_k,w_k) \in \cD_{\text{test}} } \left\|\bhH_{\SO}(\mu_k\alpha,\beta) - w_k\right\|^2,$$
where $\bhH_{\SO}$ is constructed using only the training data. 
However, the optimization problem~\eqref{eq:opt_problem} is non-convex, and solving it is a challenging task. Therefore, we seek to solve a relaxed problem. For this purpose, in the paper, we make a 2-D grid for the parameters $\alpha$ and $\beta$ in given intervals. Then, we seek to determine the parameters on the grid where the function $\cJ(\alpha, \beta)$ is minimized.  Nonetheless, solving the optimization problem \eqref{eq:opt_problem} needs future investigation and so we leave it as a possible future research problem.

\section{Numerical results}\label{sec:NumRes}
In this section, we illustrate the efficiency of the proposed methods via several numerical examples, arising in various applications. All the simulations are done on a CPU 2.6 GHz \intel~\coreifive, 8 GB 1600 MHz DDR3, \matlab~9.1.0.441655 (R2016b).
%
%
\subsection{Demo example}
At first, we discuss an artificial example to illustrate the proposed method. Let us consider a $\SO$~system of order $n=2$, $\bSigma_{\SO} = (\bM,\bD,\bK,\bB,\bC)$ whose matrices are given by:
\begin{align*}
\bM = \begin{bmatrix}
1 & 0 
\\ 
0 & 1
\end{bmatrix}, \quad \bK = \begin{bmatrix}
1 & 0 \\ 0 & 2
\end{bmatrix}, \quad \bD = \alpha\bM +\beta\bM, \quad \text{and} ~\,   \bB^T = \bC = \begin{bmatrix}
2 & 3
\end{bmatrix},
\end{align*}
with $\alpha   = 0.01$ and $\beta    = 0.02$. We collect $20$ samples $(\sigma_j, \bhH_{\SO}(\sigma_j))$, for $\sigma_j \in \iota[10^{-1},10^1]$ logarithmically spaced.  Then, we construct the \FO~and \SO~Loewner matrices in \eqref{eq:ClassicLoew} and \eqref{eq:SOLoew}, receptively.  

In Figure \ref{fig:Simple_SVD}, we plot the decay of the singular values of the $\Loew$ and $\Loew^{\SO}$ matrices.
 It can be observed that $\rank{\Loew} = 4$ and $\rank{\Loew^{\SO}} = 2$, as expected. Indeed, the demo system has a minimal $\SO$ realization of order 2 and a minimal $\FO$ realization of order 4. By applying the SVD procedure, we construct two reduced-order models of order 2, one for  $\FO$ and the other for $\SO$.  We compare the transfer functions of the original and reduced-order systems, and the results are plotted in Figure \ref{fig:Simple_Bode_Err}. The figure shows that the error between the original and $\SO$ reduced-order system is of the level of machine  precision, which means that the $\SO$ approach has recovered an equivalent realization of the original model. Additionally, the $\FO$ reduced system of order 2 was not able to mimic the same behavior of the original system, showing that a larger order is required in this case.  

\begin{figure}[tb]
	\centering
	\begin{tikzpicture}
	\begin{customlegend}[legend columns=-1, legend style={/tikz/every even column/.append style={column sep=1cm}} , legend entries={Original system, \FO~Loewner, \SO~Loewner }, ]
	\addlegendimage{black, line width = 2.0pt}
	\addlegendimage{red,dashed, line width=2.0pt}
	\addlegendimage{blue,dashed, line width=2.0pt}
	\end{customlegend}
	\end{tikzpicture}
	\centering
	\newlength\fwidth
	\newlength\fheight
	\setlength\fheight{2cm}
	\setlength\fwidth{.5\textwidth}
%
%
\begin{tikzpicture}

\begin{axis}[%
width=0.951\fwidth,
height=\fheight,
at={(0\fwidth,0\fheight)},
scale only axis,
xmin=1,
xmax=10,
xmajorgrids,
ymode=log,
ymin=1e-17,
ymax=1,
yminorticks=true,
ymajorgrids,
yminorgrids,
axis background/.style={fill=white},
title style={font=\bfseries},
ylabel={Relative singular values},
legend style={legend cell align=left, align=left, draw=white!15!black}
]
\addplot[color=red,mark=*, mark repeat = 1, line width = 1.2pt,mark size = 1.5, dashed]
  table[row sep=crcr]{%
1	1\\
2	0.743143043963565\\
3	0.322885746776308\\
4	0.301324765576857\\
5	2.56047455589109e-16\\
6	1.80910051475452e-16\\
7	6.72119359972139e-17\\
8	6.33014365017304e-17\\
9	5.10166025047229e-17\\
10	3.60565261737236e-17\\
};

\addplot [color=blue,mark=*, mark repeat = 1, line width = 1.2pt,mark size = 1.5, dashed]
  table[row sep=crcr]{%
1	1\\
2	0.526619381764003\\
3	2.59355156096016e-15\\
4	3.11297714401297e-16\\
5	1.66323635682343e-16\\
6	1.12467577712158e-16\\
7	6.14539460335997e-17\\
8	4.82071792536231e-17\\
9	1.77393396012394e-17\\
10	9.41264342927008e-18\\
};

\end{axis}
\end{tikzpicture}%
	\caption{Demo example: Decay of the singular values for the \FO~and \SO~Loewner matrices.}
	\label{fig:Simple_SVD}
	\centering
	\setlength\fheight{3cm}
	\setlength\fwidth{0.40\textwidth}
	\input{figures/Simple_Bode.tex}
	\input{figures/Simple_Err.tex}
	\caption{Demo example: The figure on the left shows the Bode plot of the original system and the $\FO$ and $\SO$ reduced-order models. The figure on the right shows  the Bode plot of the error between the original and reduced-order systems.}
	\label{fig:Simple_Bode_Err}
\end{figure}

\subsection{Building example}
Let us now consider the building model from the SLICOT library \cite{morChaV02}. It describes the displacement of a multi-storey building, for example, during an earthquake. It is a \FO~system of order $r = 48$, whose dynamics comes from a mechanical system. The Rayleigh damping coefficients here are $\alpha \approx 0.4947$ and $\beta \approx 0.0011$. 

For this example, we collect 200 samples $\bH(i\omega)$, with $\omega\in [10^0, 10^2]$. Then, we build the \FO~and \SO~Loewner matrices in \eqref{eq:ClassicLoew} and \eqref{eq:SOLoew}, receptively.  Additionally, using the heuristic procedure in Subsection \ref{subsec:OptPar}, we constructed the reduced model assuming we do not know a priori the parameters $\alpha$ and $\beta$. After this procedure, we obtain $\alpha^* = 0.495$ and $\beta^* =  0.001$, which are fairly close to the original parameters.  

In Figure \ref{fig:Build_SVD}, we plot the decay of the singular values of the \FO~Loewner matrix, the \SO~Loewner matrix for the original parameters $\alpha$ and $\beta$, and  the \SO~Loewner matrix for the estimated parameters $\alpha$ and $\beta$. The decay of the singular values for the \SO~Loewner matrix with original parameters is  faster than for the \FO~Loewner matrix. However, for the  \SO~Loewner matrix with estimated parameters, the decay of singular values starts fast and then becomes slower. This shows that if the parameters $\alpha$ and $\beta$ are not well identified, a higher reduced-order will be needed to interpolate the data.
By applying the SVD procedure, we construct three reduced-order models of order~16.  We compare the transfer functions of the original and reduced-order systems, and the results are plotted in Figure~\ref{fig:Build_Bode_Err}. This figure shows that for the \SO~Loewner approach (original parameters or with estimated parameters) outperform the classical Loewner framework. 
\begin{figure}[!tb]
	\centering
	\begin{tikzpicture}
	\begin{customlegend}[legend columns=2, legend style={/tikz/every even column/.append style={column sep=1cm}} , legend entries={Original system, \FO~Loewner, \SO~Loewner, \SO~Loewner opt. par. }, ]
	\addlegendimage{black, line width = 2.0pt}
	\addlegendimage{red,dashed, line width=2.0pt}
	\addlegendimage{blue,dashed, line width=2.0pt}
	\addlegendimage{olive,dashed, line width=2.0pt}
	\end{customlegend}
	\end{tikzpicture}
	\centering
	\setlength\fheight{3cm}
	\setlength\fwidth{.5\textwidth}
%
%
\begin{tikzpicture}

\begin{axis}[%
width=0.951\fwidth,
height=\fheight,
at={(0\fwidth,0\fheight)},
scale only axis,
xmin=0,
xmax=80,
ymode=log,
ymin=1e-15,
ymax=1,
yminorticks=true,
xmajorgrids,
ymajorgrids,
yminorgrids,
ylabel style={font=\color{white!15!black}},
ylabel={Relative singular values},
axis background/.style={fill=white},
legend style={legend cell align=left, align=left, draw=white!15!black}
]
\addplot [color=red,mark=*, mark repeat = 5, line width = 1.2pt,mark size = 1.5, dashed]
  table[row sep=crcr]{%
1	1\\
2	0.884262903620904\\
3	0.198533265599783\\
4	0.189039352244615\\
5	0.114930455468877\\
6	0.10517301370852\\
7	0.0468632968658849\\
8	0.0418702068520876\\
9	0.0215441264598517\\
10	0.0215131165032607\\
11	0.0126760413190216\\
12	0.0116297785532624\\
13	0.00640104795390445\\
14	0.00585950643622486\\
15	0.00331828387800651\\
16	0.00322466304328082\\
17	0.002359185282412\\
18	0.00207874371495311\\
19	0.00085539441108782\\
20	0.000762110914894612\\
21	0.000529266453586927\\
22	0.000496986223421157\\
23	0.000397950730131928\\
24	0.000369417224344628\\
25	0.000105765810154231\\
26	9.00258362047486e-05\\
27	4.08353871086424e-05\\
28	3.99264104160069e-05\\
29	2.06499365304193e-05\\
30	1.91816571636322e-05\\
31	7.61683794996352e-06\\
32	7.58642852146733e-06\\
33	5.19721353283549e-06\\
34	5.03548501322767e-06\\
35	3.29693851627516e-06\\
36	2.80783119053536e-06\\
37	1.99059159082581e-06\\
38	1.97201067260538e-06\\
39	1.58672048394582e-06\\
40	1.31189629948975e-06\\
41	1.44550669435192e-07\\
42	1.41172154836955e-07\\
43	1.02041139889076e-07\\
44	9.66434767259351e-08\\
45	3.96595931661071e-08\\
46	3.83442701430726e-08\\
47	1.77465622982493e-08\\
48	1.59227510091905e-08\\
49	8.66198022789956e-15\\
50	8.58170663548834e-15\\
51	8.05592218394424e-15\\
52	7.99303647902048e-15\\
53	7.15464486777106e-15\\
54	6.93241547883027e-15\\
55	6.37849353463566e-15\\
56	6.30424612746511e-15\\
57	5.99396980380498e-15\\
58	5.94299451617801e-15\\
59	5.19098522656581e-15\\
60	5.07276099378884e-15\\
61	3.87589893073772e-15\\
62	3.52648703891757e-15\\
63	3.52413840897449e-15\\
64	3.43383370472149e-15\\
65	3.41023037290115e-15\\
66	3.32419272666572e-15\\
67	3.27948458375644e-15\\
68	3.0568334073461e-15\\
69	3.01075633140084e-15\\
70	2.70257342834901e-15\\
71	2.57657753485727e-15\\
72	2.54145675347097e-15\\
73	2.49834944673205e-15\\
74	2.4863298628642e-15\\
75	2.1135141253778e-15\\
76	2.06876345776519e-15\\
77	2.02962180798374e-15\\
78	2.00232111192499e-15\\
79	1.99044296680294e-15\\
80	1.93948420468967e-15\\
};

\addplot [color=blue,mark=*, mark repeat = 5, line width = 1.2pt,mark size = 1.5, dashed]
  table[row sep=crcr]{%
1	1\\
2	0.22090689029008\\
3	0.0469236753303077\\
4	0.0337310410942981\\
5	0.0105301063566017\\
6	0.00288828577521322\\
7	0.00197352437939834\\
8	0.000792465968531127\\
9	0.000397096021312053\\
10	0.000160585847264434\\
11	9.35487629099284e-05\\
12	7.08602339754908e-05\\
13	1.54835993820348e-05\\
14	1.18448723454163e-05\\
15	5.74834998245739e-06\\
16	3.25937499620503e-06\\
17	1.70261699545365e-06\\
18	8.62554196746859e-07\\
19	4.79010028629945e-07\\
20	2.39313688495068e-07\\
21	5.36462516402034e-08\\
22	3.30344475354883e-08\\
23	1.35089086920673e-08\\
24	1.19883794913288e-08\\
25	1.00332993079932e-08\\
26	4.84721811848592e-09\\
27	3.46996123154047e-09\\
28	2.73941513575543e-09\\
29	2.04604789559544e-09\\
30	1.54041448195209e-09\\
31	1.149087368561e-09\\
32	9.67132696979096e-10\\
33	7.79292056986428e-10\\
34	5.88753920083642e-10\\
35	4.40261772718541e-10\\
36	3.51461892863733e-10\\
37	2.90505130858863e-10\\
38	2.25209693587488e-10\\
39	1.70610115902223e-10\\
40	1.62068705204421e-10\\
41	1.32147846047326e-10\\
42	1.21658908042969e-10\\
43	1.08610006098363e-10\\
44	9.10720653038063e-11\\
45	8.47152682159928e-11\\
46	7.5097293329035e-11\\
47	6.54833400256353e-11\\
48	6.43772593630632e-11\\
49	5.02112122626173e-11\\
50	4.73630749265087e-11\\
51	4.07588646621277e-11\\
52	3.34116602202412e-11\\
53	3.22007547662036e-11\\
54	2.45807871132857e-11\\
55	2.40007282568238e-11\\
56	2.09479334594847e-11\\
57	1.90921704969142e-11\\
58	1.86469059060348e-11\\
59	1.63601119628985e-11\\
60	1.53295796314996e-11\\
61	1.21644501061476e-11\\
62	1.16643977758555e-11\\
63	9.34789682750469e-12\\
64	8.54824950542895e-12\\
65	8.21147910150439e-12\\
66	7.65042162524088e-12\\
67	7.24021886211058e-12\\
68	5.82550763734518e-12\\
69	5.78109564001273e-12\\
70	5.73927255323122e-12\\
71	5.52706826101893e-12\\
72	4.59213303439995e-12\\
73	4.11788309011735e-12\\
74	3.54984202774577e-12\\
75	3.29742655611156e-12\\
76	2.80404278181948e-12\\
77	2.74956039439971e-12\\
78	2.18066371038172e-12\\
79	2.02370398797478e-12\\
80	1.89484920848785e-12\\
};

\addplot [color=olive,mark=*, mark repeat = 5, line width = 1.2pt,mark size = 1.5, dashed]
  table[row sep=crcr]{%
1	1\\
2	0.220851042643679\\
3	0.0469312206517918\\
4	0.0337294250787285\\
5	0.0105373498294873\\
6	0.00288924999195211\\
7	0.00197437570866133\\
8	0.000793370273983806\\
9	0.000397225096568633\\
10	0.000160774993600896\\
11	0.000104868639310522\\
12	9.36165302390644e-05\\
13	7.0996301806859e-05\\
14	4.27182796627706e-05\\
15	2.67450155451205e-05\\
16	1.88430827706092e-05\\
17	1.54973480961721e-05\\
18	1.45263411905236e-05\\
19	1.19061511782194e-05\\
20	1.13088259950047e-05\\
21	8.74222571754585e-06\\
22	6.99539307092082e-06\\
23	6.43704086848112e-06\\
24	5.78909423658622e-06\\
25	5.46698922529531e-06\\
26	4.35686467073872e-06\\
27	3.4344589537123e-06\\
28	3.2736924327622e-06\\
29	2.89088072400257e-06\\
30	2.7653539288827e-06\\
31	2.2592190869668e-06\\
32	1.93118837451136e-06\\
33	1.83730160596155e-06\\
34	1.70446365904945e-06\\
35	1.56642550304196e-06\\
36	1.43312514915329e-06\\
37	1.25907198844436e-06\\
38	1.11445284504049e-06\\
39	1.05387364536578e-06\\
40	8.70015972605452e-07\\
41	8.41926874919789e-07\\
42	8.38671212390434e-07\\
43	7.04082309321531e-07\\
44	7.02544461338454e-07\\
45	6.50954613174436e-07\\
46	5.67780820898558e-07\\
47	5.32732861316717e-07\\
48	4.83863426501843e-07\\
49	4.76483780658483e-07\\
50	4.43312092542665e-07\\
51	4.05888266068675e-07\\
52	3.8774369702016e-07\\
53	3.26518988095962e-07\\
54	3.24029161194859e-07\\
55	2.84521646186069e-07\\
56	2.70199961165874e-07\\
57	2.64793511599282e-07\\
58	2.53071989262207e-07\\
59	2.39652825522192e-07\\
60	2.21691390787472e-07\\
61	2.1513518640221e-07\\
62	1.91455608909103e-07\\
63	1.81712735077658e-07\\
64	1.67219369343388e-07\\
65	1.52475946929269e-07\\
66	1.49212153762848e-07\\
67	1.34490643545494e-07\\
68	1.250359101429e-07\\
69	1.20815474759165e-07\\
70	1.08504695822693e-07\\
71	1.04609622975977e-07\\
72	9.65479071662881e-08\\
73	8.64870804509559e-08\\
74	8.60022937000798e-08\\
75	7.75650558391176e-08\\
76	7.19392422754065e-08\\
77	7.0331201489267e-08\\
78	6.37875455530964e-08\\
79	5.93470985447591e-08\\
80	5.67589608762944e-08\\
};

\end{axis}
\end{tikzpicture}%
	\caption{Build example: Decay of the singular values for the \FO~ Loewner matrix and for \SO~Loewner matrices.}
	\label{fig:Build_SVD}
	\centering
	\setlength\fheight{3cm}
	\setlength\fwidth{0.40\textwidth}
	\input{figures/Build_Bode.tex}
	\input{figures/Build_Err.tex}
	\caption{Build example: The figure on the left shows the Bode plot of the original system and the $\FO$ and $\SO$ reduced-order models. The figure on the right shows  the Bode plot of the error between the original and reduced-order systems.}
	\label{fig:Build_Bode_Err}
\end{figure}

\subsection{Artificial Fishtail}
As the last example, we consider the artificial fishtail model presented in \cite{saak2019comparison}. This model comes from a finite-element discretization of the continuous mechanics model of an artificial fishtail.  After discretization, the finite-dimensional system has a $\SO$ realization of order $779,232$. For this model, the Rayleigh damping is chosen with parameters $\alpha = 1.0\cdot 10^{-4}$, $\beta =  2\cdot 10^{-4}$. It is a MIMO system, but for the numerical application, here we consider only the first transfer function, i.e., from $u_1$ to $y_1$.  

For this example, we collect 200 samples $\bH(i\omega)$, with $\omega\in [10^1, 10^4]$. Then, we build \FO~and \SO~Loewner matrices in \eqref{eq:ClassicLoew} and \eqref{eq:SOLoew}, receptively. Additionally, we also compute the reduced model using the heuristic procedure in Subsection \ref{subsec:OptPar}, for which we obtain the estimated parameters $\alpha^* \approx 1.19\cdot 10^{-4}$ and $\beta^* \approx 2\cdot 10^{-4}$.

In Figure \ref{fig:Build_SVD}, we plot the decay of the singular values of the \FO~Loewner matrix, the \SO~Loewner matrix for the original parameters $\alpha$ and $\beta$ , and  the \SO~Loewner matrix for the estimated parameters $\alpha^*$ and $\beta^*$.  By applying the SVD procedure, we construct three reduced-order models of order 8.  We compare the transfer functions of the original and reduced-order systems, and the results are plotted in Figure \ref{fig:Build_Bode_Err}.  This figure shows that the \SO~Loewner approach with original parameters and \SO~Loewner with estimated parameters outperform the classical Loewner framework.

\begin{figure}[tb]
	\centering
	\begin{tikzpicture}
	\begin{customlegend}[legend columns=2, legend style={/tikz/every even column/.append style={column sep=1cm}} , legend entries={Original system, \FO~Loewner, \SO~Loewner, \SO~Loewner opt. par. }, ]
	\addlegendimage{black, line width = 2.0pt}
	\addlegendimage{red,dashed, line width=2.0pt}
	\addlegendimage{blue,dashed, line width=2.0pt}
	\addlegendimage{olive,dashed, line width=2.0pt}
	\end{customlegend}
	\end{tikzpicture}
	\centering
	\setlength\fheight{3cm}
	\setlength\fwidth{.5\textwidth}
%
%
\begin{tikzpicture}

\begin{axis}[%
width=0.951\fwidth,
height=\fheight,
at={(0\fwidth,0\fheight)},
scale only axis,
xmin=0,
xmax=80,
ymode=log,
ymin=1e-10,
ymax=1,
yminorticks=true,
xmajorgrids,
ymajorgrids,
yminorgrids,
ylabel style={font=\color{white!15!black}},
ylabel={Relative singular values},
axis background/.style={fill=white},
legend style={legend cell align=left, align=left, draw=white!15!black}
]
\addplot  [color=red,mark=*, mark repeat = 5, line width = 1.2pt,mark size = 1.5, dashed]
  table[row sep=crcr]{%
1	1\\
2	0.973032023513909\\
3	0.00661359588794075\\
4	0.00171070984961247\\
5	0.00150592252688842\\
6	0.000286337366633571\\
7	0.000236059856873898\\
8	2.24465559784863e-05\\
9	1.39172181245341e-05\\
10	2.47180890332401e-06\\
11	2.42432847970681e-06\\
12	7.30047056258646e-07\\
13	7.11493248402737e-07\\
14	3.64631669935828e-07\\
15	3.31724370491803e-07\\
16	2.71454887854386e-07\\
17	1.35536132946986e-07\\
18	1.34216417663853e-07\\
19	1.33231280656688e-07\\
20	1.32144717958155e-07\\
21	1.1784904219648e-07\\
22	1.0279835575517e-07\\
23	1.0222748306997e-07\\
24	1.00388315173488e-07\\
25	9.94475185915504e-08\\
26	9.91801154330137e-08\\
27	9.87475383255203e-08\\
28	9.71582965225512e-08\\
29	9.69674551461604e-08\\
30	9.581104213101e-08\\
31	9.48438992048891e-08\\
32	9.31924581887597e-08\\
33	9.21754913188904e-08\\
34	8.94448769975194e-08\\
35	8.73921566976006e-08\\
36	6.68710490439495e-08\\
37	6.50307371703435e-08\\
38	6.14098631728346e-08\\
39	6.10416695417421e-08\\
40	5.70565718402732e-08\\
41	5.68023489361668e-08\\
42	5.49607382634913e-08\\
43	5.45413537878034e-08\\
44	5.1907868868802e-08\\
45	5.1907143655889e-08\\
46	5.14090769564209e-08\\
47	4.96829855775971e-08\\
48	4.86127214842128e-08\\
49	4.82739134374389e-08\\
50	4.68336813188054e-08\\
51	4.67717787147464e-08\\
52	4.47656500079581e-08\\
53	4.443078545825e-08\\
54	4.34618243017225e-08\\
55	4.271633684604e-08\\
56	4.15773867223179e-08\\
57	4.15719810977604e-08\\
58	4.06464083590899e-08\\
59	3.98999933487365e-08\\
60	3.78588508997387e-08\\
61	3.19839505083796e-08\\
62	3.10363198239521e-08\\
63	3.09820221227324e-08\\
64	3.04477681170932e-08\\
65	3.01956139003623e-08\\
66	2.98897487796287e-08\\
67	2.9295570150727e-08\\
68	2.85467689612276e-08\\
69	2.82840457992803e-08\\
70	2.32253346563092e-08\\
71	2.03867619171414e-08\\
72	2.02892891846627e-08\\
73	1.88111634056366e-08\\
74	1.80546968854368e-08\\
75	1.76067083939883e-08\\
76	1.48461429804816e-08\\
77	1.41574659821019e-08\\
78	1.26740718538602e-08\\
79	1.17271641090224e-08\\
80	1.00151805577052e-08\\
};

\addplot  [color=blue,mark=*, mark repeat = 5, line width = 1.2pt,mark size = 1.5, dashed]
  table[row sep=crcr]{%
1	1\\
2	0.000392301140911182\\
3	4.09961345224768e-05\\
4	2.70454145152963e-06\\
5	2.36632360856008e-06\\
6	1.01251257605257e-06\\
7	9.89163080936928e-07\\
8	9.32237087179955e-07\\
9	6.61927289427861e-07\\
10	5.58637051911708e-07\\
11	5.01964254084446e-07\\
12	4.58901236672335e-07\\
13	4.38223705285713e-07\\
14	3.73062828364701e-07\\
15	3.71289156589115e-07\\
16	3.47789845316191e-07\\
17	3.32403386969812e-07\\
18	2.71433286031889e-07\\
19	2.54857397590716e-07\\
20	2.3189761614394e-07\\
21	1.87113201677243e-07\\
22	1.76049119781656e-07\\
23	1.72907373672184e-07\\
24	1.54553925352405e-07\\
25	1.39912540977287e-07\\
26	1.30266006653086e-07\\
27	1.27777468122396e-07\\
28	1.17897756741257e-07\\
29	1.14819331060836e-07\\
30	1.11055102887115e-07\\
31	1.06913515991174e-07\\
32	9.80953152338396e-08\\
33	9.37635433521694e-08\\
34	8.85595262073022e-08\\
35	6.90468067936031e-08\\
36	5.78388434985905e-08\\
37	5.45160592682666e-08\\
38	5.00535507539309e-08\\
39	4.85006679615501e-08\\
40	4.53064587000199e-08\\
41	4.38089906722275e-08\\
42	3.68145691807414e-08\\
43	3.13047953202828e-08\\
44	2.42485893235069e-08\\
45	2.31477372254823e-08\\
46	2.12218738889015e-08\\
47	1.79840673905407e-08\\
48	1.53179259457375e-08\\
49	1.46645826409121e-08\\
50	1.29836160688268e-08\\
51	1.13719546060129e-08\\
52	9.9427629866088e-09\\
53	7.93955342294402e-09\\
54	6.11332249392785e-09\\
55	5.71573228431384e-09\\
56	4.67374861809876e-09\\
57	4.46689898031948e-09\\
58	4.21252135495158e-09\\
59	3.91222269181454e-09\\
60	3.7134511838012e-09\\
61	3.50567046833967e-09\\
62	3.2268191802379e-09\\
63	3.13983571397962e-09\\
64	2.89222512329638e-09\\
65	2.55786798283421e-09\\
66	2.30580414224341e-09\\
67	2.11545465603257e-09\\
68	2.01419422441589e-09\\
69	1.96001148828218e-09\\
70	1.81176377824801e-09\\
71	1.71213566285589e-09\\
72	1.64430273888682e-09\\
73	1.61224900291891e-09\\
74	1.28120930630419e-09\\
75	1.22648809204026e-09\\
76	1.19430788917116e-09\\
77	1.15356063571406e-09\\
78	1.06140139352974e-09\\
79	1.03823578327574e-09\\
80	9.14794611911553e-10\\
};

\addplot  [color=olive,mark=*, mark repeat = 5, line width = 1.2pt,mark size = 1.5, dashed]
  table[row sep=crcr]{%
1	1\\
2	0.000392296361075857\\
3	4.10007789821021e-05\\
4	2.7114810044086e-06\\
5	2.53253877870913e-06\\
6	1.00668639491441e-06\\
7	9.62832299570512e-07\\
8	9.33499055456383e-07\\
9	8.43343665033271e-07\\
10	6.5145651711267e-07\\
11	5.47705237029639e-07\\
12	4.59775340580041e-07\\
13	4.43967045680767e-07\\
14	4.2891416679469e-07\\
15	3.75301260505111e-07\\
16	3.51153733437834e-07\\
17	3.3672389006934e-07\\
18	2.5889702190386e-07\\
19	2.36295230467541e-07\\
20	2.09109495987602e-07\\
21	1.8581198884468e-07\\
22	1.786049813569e-07\\
23	1.75344205992726e-07\\
24	1.67495653213629e-07\\
25	1.3914772258521e-07\\
26	1.3123075499015e-07\\
27	1.27862904556221e-07\\
28	1.15562978554609e-07\\
29	1.14866503330025e-07\\
30	1.09956329619198e-07\\
31	1.02925690117082e-07\\
32	1.01817508467541e-07\\
33	9.22807129470714e-08\\
34	7.3836037660835e-08\\
35	7.09904052592404e-08\\
36	6.37471594607169e-08\\
37	5.66919097533171e-08\\
38	5.47825855940823e-08\\
39	5.30571614410334e-08\\
40	4.95637468477412e-08\\
41	4.34679347990894e-08\\
42	3.97607337128458e-08\\
43	3.55964769042313e-08\\
44	3.13069929891122e-08\\
45	2.97486578931728e-08\\
46	2.68824317937107e-08\\
47	2.49664459519913e-08\\
48	2.3273766780646e-08\\
49	1.87474602606276e-08\\
50	1.51710689999692e-08\\
51	1.43316665027398e-08\\
52	1.2569125643892e-08\\
53	1.13860144456154e-08\\
54	1.1335741563966e-08\\
55	9.92711720022657e-09\\
56	9.50599755395804e-09\\
57	8.96997621539502e-09\\
58	8.72841017575751e-09\\
59	8.30335344888338e-09\\
60	6.80967622043632e-09\\
61	6.76176103395216e-09\\
62	4.9921062473934e-09\\
63	4.78044949322637e-09\\
64	4.03700773439504e-09\\
65	3.78806633505497e-09\\
66	3.55781834693208e-09\\
67	3.40572079752929e-09\\
68	3.06113692273392e-09\\
69	3.00663065808601e-09\\
70	2.9174755446251e-09\\
71	2.76162714587687e-09\\
72	2.58726242144095e-09\\
73	2.45472609516107e-09\\
74	2.23726239894044e-09\\
75	1.96721898115351e-09\\
76	1.88671972170705e-09\\
77	1.69807835796999e-09\\
78	1.63168395443729e-09\\
79	1.58744207463981e-09\\
80	1.33282056144577e-09\\
};

\end{axis}
\end{tikzpicture}%
	\caption{Fishtail example: Decay of the singular values for the \FO~ Loewner matrix and for \SO~Loewner matrices.}
	\label{fig:Fish_SVD}
	\centering
	\setlength\fheight{3cm}
	\setlength\fwidth{0.40\textwidth}
	\input{figures/Fish_Bode.tex}
	\input{figures/Fish_Err.tex}
	\caption{Fishtail example: The figure on the left shows the Bode plot of the original system and the $\FO$ and $\SO$ reduced-order models. The figure on the right shows  the Bode plot of the error between the original and reduced-order systems.}
	\label{fig:Fish_Bode_Err}
\end{figure}

\section{Conclusions}\label{sec:Conc}
In this paper, we have studied the problem of the identification of Rayleigh-damped second-order systems from frequency data. To that aim, we propose modified  \SO~Loewner matrices which are the key tools to construct a realization interpolating the given data. Additionally, in the case of redundant data, an SVD-based scheme is presented to construct reduced-order models. Moreover, a heuristic optimization problem is sketched to estimate the damping parameters. Finally, we have illustrated the efficiency of the proposed approach in some numerical examples, and we compared the results with the classical Loewner framework. 
\vspace{-0.7cm}
\section*{Acknowledgement}
This work was supported by \emph{Deutsche Forschungsgemeinschaft (DFG)}, Collaborative Research Center CRC 96 "Thermo-energetic Design of Machine Tools".

\vspace{-0.5cm}
\bibliographystyle{siam} 
\bibliography{igorBiblio,mor,csc}
\end{document}